\documentclass[11pt, letterpaper]{article}    % Enable this line and disable the
\usepackage{graphicx,setspace,fullpage}
\usepackage[tbtags]{amsmath}                % added by editor
\usepackage{amssymb}                % added by editor
\usepackage{latexsym}               % added by editor
\newcommand{\pfbox}{\hfill\mbox{$\Box$}}

\newcommand{\R}{\mbox{$\mathbb{R}$}}

\newtheorem{lemma}{Lemma}
\newtheorem{remark}{Remark}
\newtheorem{theorem}{Theorem}

\newtheorem{cor}{Corollary}

\begin{document}
%\setstretch{.88}
%\input{BlocksStyle}
%\input{SumNodesStyle}
%\input{HighlightingStyle}
%\input{LinesStyle}

%\begin{frontmatter}

\title{Consensus Control for Heterogeneous Multi-Agent Systems}%\thanksref{footnoteinfo}}
\author{Luis D. Alvergue, Abhishek Pandey, Guoxiang Gu, and Xiang Chen}
\date{May 6, 2015}
\maketitle
% Title, preferably not more

%\thanks[footnoteinfo]{Corresponding author L. Alvergue. This
%research is supported in part by the US Air Force and NASA.}

%\author[LAAPGG]{Luis D. Alvergue}\ead{lalver1@tigers.lsu.edu},    % Add the
%\author[LAAPGG]{Abhishek Pandey}\ead{apande3@tigers.lsu.edu},               % e-mail address
%\author[LAAPGG]{Guoxiang Gu}\ead{ggu@lsu.edu},  % (ead) as shown
%\author[XC]{Xiang Chen}\ead{xchen@uwindsor.ca}  % (ead) as shown

%\address[LAAPGG]{%Division of Electrical and Computer Engineering,
%School of Electrical Engineering and Computer Science, Louisiana
%State University, Baton Rouge, LA 70803, USA.}  % Please supply
%\address[XC]{Department of Electrical and Computer Engineering, University of Windsor, Windsor, ON, N9B 3P4, Canada.}             % full addresses
%\address[GG]{The White House, Baiae}        % here.

%\begin{keyword}                           % Five to ten keywords,
%Distributed control, multi-agent systems,
%and output consensus                    % chosen from the IFAC
%\end{keyword}                             % keyword list or with the
                                          % help of the Automatica
                                          % keyword wizard

\begin{abstract}                          % Abstract of not more than 200 words.
We study distributed output feedback control for a
heterogeneous multi-agent system (MAS), consisting of
$N$ different continuous-time linear dynamical systems.
For achieving output consensus, a virtual
reference model is assumed to generate the desired
trajectory that the MAS is required to track
and synchronize. A distinct feature of our results
lies in the local optimality and robustness achieved by our
proposed consensus control algorithm.
In addition our study is focused
on the case when the available output measurements
contain only relative information from
the neighboring agents and reference signal.
Indeed by exploiting properties of strictly
positive real (SPR) transfer
matrices, conditions are derived
for the existence of distributed
output feedback control protocols, and
solutions are proposed to synthesize the stabilizing and
consensus control protocol over a given
connected digraph. It is shown that design
techniques based on the LQG,
LQG/LTR and ${\mathcal H}_{\infty}$ loop shaping
can all be directly applied to synthesize
the consensus output feedback control protocol,
thereby ensuring the local optimality and stability robustness.
Finally the reference trajectory is required to be
transmitted to only one or a few agents
and no local reference models are employed
in the feedback controllers thereby
eliminating synchronization of the local reference models.
%Both significantly lower the communication
%overhead. In addition,
Our results complement the existing ones,
and are illustrated by a numerical example.
\end{abstract}

%\end{frontmatter}

\section{Introduction}

\noindent In the last decade, research on cooperative
MASs has intensified mainly due to the wide variety
of applications that make use of the MAS framework,
cf. \cite{jlm,lwcv,ofm,rba}. % and references therein.
Early work focused on agents with scalar integrator
dynamics. This includes \cite{jlm,rb} where the
agents' dynamics are represented by linear switched systems
or by changing interaction topologies, and
\cite{m,om}, where the dynamics correspond to
time-varying systems. Scenarios that are more
realistic, such as systems with communication
noises, have also been considered \cite{lz,rbk}.
More recently the majority of the existing work has been
concerned with homogeneous agents (i.e. agents
with identical dynamics) represented by the state-space
model. See for example \cite{ldch,mz,ss,zld}. Such
MASs are more general and include
integrator dynamics as a special case.
Solutions to the problem of achieving consensus
through distributed control protocols are presented
in these investigations. Since the separation
principle holds for MASs of homogeneous agents,
the control protocols are mostly observer-based,
including Luenberger observers.
%In addition nonlinear dynamic agents are also
%under investigation by several authors \cite{?,?}
%for passive or dissipative systems.

Motivated by the recent developments, we
study consensus control of heterogeneous
MASs which is in general more difficult
than that of homogeneous MASs. Some recent
%controller synthesis
results include \cite{gyss,kss,l,tz,wsa} %\cite{gyss,kss,l,wsa}
and references therein. The authors of
\cite{gyss,kss,wsa} smartly
bypass the difficulty posed
by heterogeneity by embedding a homogeneous
reference model into each agent's local
controller. The state vectors of the local reference models
are synchronized prior to achieving
output consensus of the heterogeneous MAS.
In addition, the results in \cite{wsa} make use
of those in \cite{m1} and \cite{ss} to prove asymptotic synchronization
over time-varying directed graphs that satisfy
a uniform connectivity condition.
The authors of \cite{gyss} study
the same consensus control problem over
a fixed or time-invariant graph topology. However it
differs from \cite{kss,wsa} in that
it considers heterogeneous non-introspective agents.
Roughly speaking, non-introspective agents
refer to those agents that make use of only relative
information with respect to neighbors or reference inputs. The
controller design is hinged on an observer that depends
on a high gain parameter. It is important to highlight that
in \cite{gyss} the network graph is only present at the output but
not in the input of the MAS. The case when the graph is present
at the input introduces more difficulties.
%Non-introspective
%agents possess no knowledge about their own state or
%output separate from what is received via the
%network.
%However the controller design
%requires right-invertibility of
%the agent's realization
%which does not hold for
%$m$-input/$p$-output systems if $m<p$.
It is worth mentioning that
analysis results are also available for
heterogeneous MASs. For instance, \cite{l}
provides a theoretical condition, stated
in terms of ``intersection dynamics",
for synchronization that requires inclusion of the
internal modes of the root agent or reference
model. The work \cite{tz} provides
a consensusability condition in the presence of
unknown communication delays. Finally,
nonlinear heterogeneous MASs
are also studied in, for instance,
\cite{cs,hh,zhl}, where the focus is on
passive and dissipative systems. Such
nonlinear agents exclude unstable dynamical systems,
yet the results are instrumental to future work
on nonlinear MASs.

In spite of the recent developments, the
design of distributed and local control protocols
to achieve not only feedback stability
but also output consensus in tracking reference
trajectories remains a major challenge. For
instance local optimality and robustness
are not addressed thus far in the existing
literature for heterogeneous MASs.
This paper is aimed at developing a more
accessible method for consensus control
and deriving a consensusability condition
for heterogeneous MASs, taking local optimal
and robust control into consideration.
It will be shown that
similar results to the ones found in
\cite{ldch,mz,zld} for homogeneous MASs
are available for heterogeneous MASs.
It will also be shown that
existing design methods, such as
linear quadratic regulator (LQR),
linear quadratic Gaussian (LQG) and
loop transfer recovery (LTR) \cite{am}, and
${\mathcal H}_{\infty}$
loop shaping \cite{mg}, developed for
multi-input/multi-output (MIMO)
feedback control systems can be employed to synthesize
consensus controllers for heterogeneous MASs.
Consequently, the consensus method as proposed in
this paper ensures the local optimality and
stability robustness. To show that our controller
design achieves output consensus, we exploit a property of SPR transfer
matrices in a feedback connection based on which optimal
and robust control can be applied.
Another distinction of our work compared to other
investigations is that we provide a solution to the consensus
problem for the case when not all agents have
access to the reference trajectory. In fact, if the
communication graph is connected (or contains a
spanning tree), then it is sufficient for one agent to have
access to the reference output for the heterogeneous MAS
to synchronize and to achieve consensus. Moreover our proposed
design method does not require duplication of
the reference model in each of the local controllers thereby
eliminating synchronization of the local reference models
commonly adopted in the existing work.
Both features lower significantly the communication
overhead between agents by avoiding the need to
communicate the reference trajectory to all agents
and by removing additional synchronization
between local reference models. Furthermore
we focus on non-introspective agents as in \cite{gyss}
and use only relative information for both
state feedback and state estimation.
%Since the
%controller gains are computed based on either
%$\mathcal H_\infty$ loop shaping or LGQ/LTR methods,
%each controlled agent is robust to perturbations
%in the form of coprime factor %uncertainties or
%or gain/phase uncertainties, respectively.
%In addition, in contrast to
%much of the current work, with the exception of \cite{tm,yt},
%our approach relies heavily upon both transfer function
%and state space description of the collective MAS dynamics.
%The output consensus control is decomposed into
%the one for state feedback and the
%other for state estimation.
%By applying a limiting argument to a
%specific transfer function, we are able
%to state conditions for output consensus to a given
%reference trajectory.

The notation in this paper is more or
less standard. The $N$-dimensional real/complex space
is denoted by $\mathbb{R}^N/\mathbb{C}^N$. % with
%${\rm 1}_N\in\mathbb{R}^N$ for a vector having
%all its entries equal to one.
The space of all $p\times m$ real/complex matrices
is denoted by $\mathbb{R}^{p\times m}/\mathbb{C}^{p\times m}$.
Let $M=\left[\ \mu_{ij}\ \right]$ be a matrix with $\mu_{ij}$
the $(i,j)$th entry. Its $i^{th}$ singular value
is denoted by $\sigma_i(M)$ arranged in
descending order with $\overline\sigma(M)=\sigma_1(M)$.
For square $M$, its $i^{th}$ eigenvalue is denoted by
$\lambda_i(M)$. The symbol $\otimes$ represents the
Kronecker product.
A real square matrix $M$ is called {\it row dominant} if
$|\mu_{ii}|\geq\sum\limits_{j\neq i}|\mu_{ij}|$,
{\it column dominant} if $|\mu_{jj}|\geq\sum\limits_{i\neq j}|\mu_{ij}|$,
and doubly dominant if it is both row and column dominant.
If the inequalities are strict then one calls such
matrices strictly row or column or doubly dominant.
The rest of the notation will be made clear as we proceed.

\section{Preliminaries}

This section prepares the results in later sections
by reviewing some graph theory, formulating the consensus problem,
and providing an important preliminary result.

\subsection{Graph and Its Associated Matrices}

We focus on directed graphs (digraphs), although our
results also hold for MASs over an undirected feedback graph.
Consider a weighted digraph specified by $\mathcal
G=\left(\mathcal V,\mathcal E\right)$, where
$\mathcal V=\left\{v_i\right\}_{i=1}^N$ is the set
of nodes and $\mathcal E\subset\mathcal V\times\mathcal V$
is the set of edges or arcs, where an edge starting at node $i$
and ending at node $j$ is denoted by $\left(v_i,v_j\right)\in\mathcal E$.
The node index set is denoted by $\mathcal N=\left\{1,\cdots,N\right\}$.
The neighborhood of node $i$ is denoted by the set $\mathcal N_i
=\left\{j\;|\;(v_i,v_j)\in\mathcal E\right\}$. A path on the digraph
is an ordered set of distinct nodes
$\left\{v_{i_1},\cdots,v_{i_K}\right\}$ such that
$(v_{i_{j-1}},v_{i_j})\in\mathcal E$. If there is a path in
$\mathcal G$ from node $v_i$ to node $v_j$, then $v_j$ is said
to be reachable from $v_i$, denoted as $v_i\rightarrow v_j$.
The digraph is called strongly connected if $v_i\rightarrow v_j$
and $v_j\rightarrow v_i\;\forall\;i,j\in\mathcal N$.
The set of nodes which can reach node $v_k$ is denoted as
$\mathcal S_k=\left\{v_j\in\mathcal V:\;\exists
\mbox{ a path } v_j\rightarrow v_k\right\}$.
The digraph is called connected if there exists a node $v_k$ such that
$v_j\in\mathcal S_k$ for $j=1,\ldots,N,\;j\neq k$.
%If $v_i\rightarrow v_j\;\forall\;j\in\mathcal N$, then
%$v_i$ is called a connected node.
%The digraph is called connected if there exists a connected node in $\mathcal G$.
%{\bf If a node $v_i$ is reachable from every other node in the digraph,
%then it is globally reachable. Lemma: A digraph has a globally reachable node
%if and only if 0 is a simple eigenvalue of $\mathcal L$ \cite{lfm}.}

%\subsection{Matrices Associated with Graphs}
Let $\mathcal A=\left[\ a_{ij}\ \right]\in \mathbb{R}^{N\times N}$
be the weighted adjacency matrix.
%We assume constant $\mathcal A$ in this paper.
The value of $a_{ij}\geq 0$ represents the coupling strength of edge
$(v_i,v_j)$. % and $a_{ij}>0$ if $e_{ij}\in\mathcal E$, otherwise $a_{ij}=0$.
Self edges are not allowed, i.e., $a_{ii}=0\;\forall\;i\in\mathcal N$.
Denote the degree matrix for $\mathcal A$ by $\mathcal D
={\rm diag}\left\{{\rm deg}_1,\cdots,{\rm deg}_N\right\}$ with
${\rm deg}_i=\sum\limits_{j\in\mathcal N_i}a_{ij}$ and the
Laplacian matrix as $\mathcal L=\mathcal D-\mathcal A$.
Let ${\rm 1}_N\in\mathbb{R}^N$ be a vector of 1's.
It is clear that $\mathcal L1_N=0$ and thus
it has at least one zero eigenvalue.
It is also known that ${\rm Re}\{\lambda_i(\mathcal L)\}
\geq 0\;\forall\;i$. In fact the only eigenvalues of the
Laplacian matrix on the imaginary axis are zero
in light of the Gershgorin circle theorem. In addition,
%for directed graphs,
zero is a simple eigenvalue of $\mathcal L$, if and only if
$\mathcal G$ is a connected digraph.
We would like to call attention to the fact that similar conditions on
the eigenvalues of $\mathcal L$ can be obtained through other
properties of $\mathcal G$. For example, in \cite{rb}, it is stated that
$\mathcal L$ has one zero eigenvalue if and only
if $\mathcal G$ has a spanning tree.

An $M$-matrix has all its
off-diagonal elements being either negative or zero,
and all its principal minors being strictly positive.
A semi $M$-matrix differs from the $M$-matrix in that
it has all its principal minors being
nonnegative. Clearly the Laplacian matrix
is a semi $M$-matrix. An $M$-matrix is
said to be row/column (strictly) dominant,
if each of its rows/columns sums to a (strictly)
positive number. More properties on $M$-matrices
may be found in \cite{w}.

%\begin{fact}\label{F1}
%Suppose that $M\in \R^{N\times N}$ has all its
%off-diagonal elements being either negative or zero.
%Then the following statements are equivalent:
%\begin{enumerate}
%    \item $M$ is an $M$-matrix;
%    \item $-M$ is Hurwitz;
%    \item The leading principal minors of $M$ are all positive;
%    \item There exists $D = {\rm diag}\{d_1, \cdots, d_N\}$
%    such that $MD$ (resp. $DM$) is strictly row (resp. column) dominant.
%\end{enumerate}
%\end{fact}

\subsection{Problem Formulation}

We consider $N$ heterogeneous agents
with the dynamics of the $i$th agent described by
\begin{align}\label{ss}
    \dot{x}_i(t) = A_ix_i(t)+B_iu_i(t),
    \; \; \; \; y_i(t)=C_ix_i(t)
\end{align}
where $x_i(t)\in\mathbb{R}^{n_i}$ is the state,
$u_i(t)\in\mathbb{R}^{m_i}$ is the control input,
and $y_i(t)\in\mathbb{R}^{p}$ is the controlled output.
Note that the state dimension $n_i$ and input
dimension $m_i$ can be different from each other.
However, all agents have the same number of outputs,
and we assume that $p\leq m_i$ for all $i$. Let
$I_n$ be the $n\times n$ identity matrix.
The $i$th agent admits transfer matrix
$P_i(s)=C_i(sI_{n_i}-A_i)^{-1}B_i$.

For heterogenous MASs, we are concerned with
output consensus aimed at achieving
\begin{align}%\label{yiyjc}
    \lim_{t\rightarrow \infty}
    \left[y_i(t)-y_j(t)\right] = 0 \; \; \forall \ i, j\ \in {\mathcal N}.
\label{Eq:OCobj}
\end{align}
In order to take performance into consideration,
our focus will be synchronization of all
the agents' output to the desired trajectories
generated by an exosystem or reference model described by
\begin{align}\label{rm}
    \dot{x}_0(t) = A_0x_0(t), \; \; \; \;
    y_0(t) = C_0x_0(t)
\end{align}
with zero steady-state error. This is a virtual reference
generator where $x_0(t)\in\R^{n_0}$ and $y_0(t)\in\R^p$, and all eigenvalues of $A_0$ are restricted to lie on the imaginary axis.
A real-time reference trajectory may not be actually from
this exosystem, but consists of piece-wise step, ramp,
sinusoidal signals, etc. whose poles coincide with
eigenvalues of $A_0$. To reduce the communication overhead,
the reference signal is often
transmitted to only one or a few of the $N$ agents. Following
\cite{zld}, we call these agents controlled agents.

Assume that the realizations of $N$ agents
are all stabilizable and detectable.
In this paper we will
study under what condition
for the feedback graph,
there exist distributed stabilizing controllers and
consensus control protocols such that the
outputs of $N$ agents satisfy
not only (\ref{Eq:OCobj}) but also
\begin{align}\label{yi0}
    \lim_{t\rightarrow \infty}
    \left[y_i(t)-y_0(t)\right] = 0 \; \; \forall \ i\ \in {\mathcal N}.
\end{align}
Moreover we will study how
to synthesize the required distributed and
local controllers in order to achieve output consensus,
taking performance into account.

\subsection{A Fundamental Lemma}

Let $e_i\in\mathbb{R}^N$ be a vector
with 1 in the $i$th entry and zeros elsewhere.
We state the following lemma that is instrumental
to the main results of this paper. %Essentially, it is a technical
%condition that determines the network graph signals that
%are available for the distributed control protocols that we will propose. % \cite{apgc}.

\begin{lemma}\label{le1}
Suppose that ${\mathcal L}$ is the Laplacian matrix
associated with the directed graph $($digraph$)$ ${\mathcal G}$.
The following statements are equivalent:
\begin{verse}
{\rm (i)}\ There exists an index $i_R\in{\mathcal N}$
such that
\begin{align}\label{rankcond}
    {\rm rank}\left\{{\mathcal L}+e_{i_R}e_{i_R}^{\prime}\right\}=N;
\end{align} \\
{\rm (ii)}\ There exist diagonal $D>0$ and an index $i_R\in{\mathcal N}$
such that %with $G={\rm diag}(e_{iR})$,
\begin{align}\label{LE1}
    {\mathcal M}_D + {\mathcal M}_D^{\prime}>0,
    \; \; \; {\mathcal M}_D=D({\mathcal L}+e_{iR}e_{iR}');
\end{align} \\
%\item
{\rm (iii)}\ The digraph ${\mathcal G}$ is connected.
\end{verse} %description}
\end{lemma}

Proof: Let $\Rightarrow$ stand for ``implies''. We
will show that (iii) $\Rightarrow$ (i)
$\Rightarrow$ (ii) $\Rightarrow$ (iii) in order to establish
the equivalence of the three statements. For
(iii) $\Rightarrow$ (i), assume that
${\mathcal G}$ is connected. Then there exists
a reachable node $v_{i_R}\in{\mathcal V}$
for some index $i_R\in{\mathcal N}$.
We adapt a result in \cite{hh} (Lemma 2) to
construct an augmented graph $\overline{{\mathcal G}}$
by adding a node $v_0$, and adding an edge from $v_{i_R}$
to $v_0$ with weight $1$. The augmented graph is again
connected with $v_0$ as the only reachable node.
It follows that the Laplacian matrix associated with
the augmented graph $\overline{{\mathcal G}}$ is given by
\begin{align*}
    \overline{{\mathcal L}}=\left[
    \begin{array}{cc}
    0 & \begin{array}{ccc}
    \cdots & & 0\end{array} \\
    -e_{i_R}& {\mathcal L}+e_{i_R}e_{i_R}^\prime
    \end{array}\right].
\end{align*}
%which assures the out-degree of $v_{iR}$ is not changed.
Since the augmented graph is connected,
the Laplacian matrix $\overline{{\mathcal L}}$
has only one zero eigenvalue, implying
the rank condition (\ref{rankcond}), and thus (i) is true.

For (i) $\Rightarrow$ (ii), assume that
the rank condition (\ref{rankcond}) is true. Then
${\mathcal L}+e_{iR}e_{iR}^\prime$
is an $M$-matrix, because it is not only
a semi $M$-matrix but also has all its
eigenvalues on strict right
half plane, in light of the Gershgorin circle theorem.
Using the properties of $M$-matrices in \cite{w}, we conclude
the existence of a diagonal matrix $D$ such that
\[
    {\mathcal M}_D = D({\mathcal L} + e_{i_R}e_{i_R}^\prime)
    %= D{\mathcal L}+G
\]
is strictly column dominant. % where
%$G=De_{i_R}e_{i_R}^\prime$ is diagonal and has
%rank 1.
Since ${\mathcal M}_D$ is row dominant,
although not strictly, ${\mathcal M}_D+{\mathcal M}_D^\prime$ is both
strictly row and column dominant, thereby
concluding (ii).

For (ii) $\Rightarrow$ (iii), assume that
(\ref{LE1}) is true. Then ${\mathcal M}_D$
is an $M$-matrix, by the fact that all
its eigenvalues lie on strict right half plane.
Hence there holds
\[ %begin{align*}
    N = {\rm rank}\{D^{-1}{\mathcal M}_D\}
    ={\rm rank}\{{\mathcal L}+e_{i_R}e_{i_R}^\prime\}
    \leq {\rm rank}\{{\mathcal L}\}+1,
    %N&={\rm rank}\{D^{-1}{\cal M}\} \\
    %&={\rm rank}\{{\cal L}+ge_{i_R}e_{i_R}^\prime\}\\
    %&\leq {\rm rank}\{{\cal L}\}+1,
\] %end{align*}
by the rank inequality and
${\rm rank}\{e_{i_R}e_{i_R}^\prime\}=1$.
The above implies that
${\rm rank}\{{\mathcal L}\} \geq N-1$.
It follows that
the Laplacian matrix ${\mathcal L}$ has only
one zero eigenvalue, concluding that the
graph ${\mathcal G}$ is connected,
and thus (iii) is true. The proof is now complete.
\pfbox

%Although Lemma \ref{le1}
%considers only the directed graph,
%the result holds for the
%undirected graph with a similar and
%much simpler proof. %The requirement that a digraph
%is strongly connected has been reported
%in several papers to ensure consensusability
%\cite{cpbm,jlm}, whereas
%the condition (\ref{ranki}) with the addition
%of the connected condition is new and crucial to the
%main results of this paper.

\begin{remark}\label{rem1}
{\em For any $i_R\in{\cal N}$ corresponding to
a reachable node, there exists a diagonal matrix $D>0$ such that
${\mathcal M}_D+{\mathcal M}_D'>2I$, i.e.,
\begin{align}\label{kap}
    %{\mathcal M}+{\mathcal M}' =
    D{\mathcal L}+{\cal L}'D
    +2De_{i_R}e_{i_R}^\prime > 2I
\end{align}
Efficient algorithms for linear matrix inequality (LMI)
can be used to search for $D$. Such a $D$ helps
to design control protocol achieving not only
the MAS feedback stability but also optimizing
local performance for each agent.
However the existence of the stabilizing control
protocols does not depend on $D$.
For MIMO agents with $m$-input/$p$-output,
a commonly adopted graph has the weighted adjacency matrix
in the form of ${\mathcal A} = \{a_{ij}I_{q_i}\}$
with $q_i \equiv m_i$ or $q_i \equiv p$. Thus $D$ is modified to
$D={\rm diag}(d_1I_{q_1},\cdots,d_NI_{q_N})$
and ${\cal M}_D$ has a similar modification.
%\begin{align}\label{dgd}
%    \begin{array}{rcl}
%    D&=&{\rm diag}(d_1I_q,\cdots,d_NI_q), \\
%    G&=&{\rm diag}(g_1I_q,\cdots,g_NI_q)
%    \end{array}
%%\end{align}
%with only one nonzero $g_i>0$.
}\end{remark}

\section{Distributed Stabilization}\label{Sec:diststab}

This section is focused on the distributed control
protocol over the connected graph ${\mathcal G}$,
represented by its Laplacian matrix ${\mathcal L}$.
Distributed stabilization will be studied
and a stabilizability condition will be derived
for both the case of state feedback and output feedback.

%\subsection{State Feedback}

%Distributed local state feedback assuming $C_i=I_{n_i}$.
Let $v_i(t)=F_{0i}x_0(t)-F_ix_i(t)$ be the full
information control signal for the $i$th agent
with $F_i$ the state feedback gain for
agent $i$, and $F_{0i}$ the state feedforward gain
from the reference model.
Denote $r_i(t)=F_{0i}x_0(t)$ and $\delta_{\rm K}(\cdot)$
as the Kronecker delta function.
Consider the full information (FI) control protocol
for the $i$th agent specified by
\begin{align}\label{sfu}
    u_i(t) = d_i\delta_{\rm K}(i-i_R)v_i(t)+d_i
    \sum_{j=1}^N a_{ij}\left[v_i(t)-v_j(t)\right]
    %- g_i(F_ix_i-C_0x_0) %F_{0i}x_0)
\end{align}
with $i_R$ corresponding to one of the reachable nodes
in order to minimize the communication overhead.
This is why Lemma \ref{le1} becomes useful.
By denoting $x(t)$ as the collective state and
$r(t)$ as the collective reference, i.e.,
the stacked vector of $\{x_i(t)\}_{i=1}^N$
and $\{r_i(t)\}_{i=1}^N$, respectively,
the closed loop dynamics with protocol (\ref{sfu}) can be written as
\begin{align}\label{sfclu}
    \dot{x} = \left[A
    -B{\mathcal M}_DF\right]x +
    B{\mathcal M}_Dr.
\end{align}
where $A = {\rm diag}(A_1,\cdots,A_N)$, $B = {\rm diag}(B_1,\cdots,B_N)$,
and $F={\rm diag}(F_1,\cdots,F_N)$. %, and ${\cal M}=D({\cal L}+e_{i_R}e_{i_R}')$.
The following result is concerned
with distributed stabilization
under state feedback.
%Since $C_i=I_{n_i}$, the
%output consensus control is aimed at $F_ix_{i_{\rm ss}} = C_0x_{0_{\rm ss}}$
%for $1\leq i\leq N$. A similar analysis yields
%the condition $F_iA_ix_{i_{\rm ss}} = 0$ for each $i$.

\begin{theorem}\label{th2}
Suppose that $(A_i,B_i)$ is stabilizable for
each $i\in{\mathcal N}$.
There exist stabilizing {\rm FI} control protocols in the form of
$(\ref{sfu})$ for the feedback {\rm MAS} over
the directed graph ${\mathcal G}$, if ${\mathcal G}$ is connected.
\end{theorem}

Proof: Since $\mathcal G$ is connected,
Lemma \ref{le1} and Remark \ref{rem1} imply
the existence of a reachable node with index $i_R\in{\cal N}$ and
$D={\rm diag}(d_1I_{m_1},\cdots,d_NI_{m_N})$
such that ${\cal M}_D+{\cal M}_D'>2I$ holds.
Feedback stability of the underlying MAS requires
the existence of $F$ such that
\begin{align}\label{rfd}
    %\lambda(s):=
    {\rm det}\left(sI-A+
    B{\cal M}_DF\right)\neq 0
    \; \; \;\forall \ {\rm Re}\{s\}\geq 0.
\end{align}
%Denote $\underline{D}={\rm diag}(d_1I_{n_1},\cdots,d_NI_{n_N})$.
Let $Z={\cal M}_D-I$. Then $Z+Z'>0$, and thus
the inequality (\ref{rfd}) is equivalent to
\begin{align*} %\label{eq:detstab}
    {\rm det}\left(sI-A+BF+BZF\right)\neq 0
    \; \; \; \forall \ {\rm Re}\{s\}\geq 0.
\end{align*}
Denote $T_F(s)=F(sI-A+BF)^{-1}B$.
The above inequality is in turn equivalent to
\begin{align}\label{prc}
    {\rm det}\left[I+T_F(s)Z\right]\neq 0\;
    \; \; \forall \ {\rm Re}\{s\}\geq 0.
\end{align}
To show the existence of the stabilizing
$F={\rm diag}(F_1,\cdots,F_N)$ that satisfies
inequality (\ref{prc}), consider the algebraic
Riccati equation (ARE)
\begin{align}\label{arexi}
    A_i'X_i+X_iX_i-XB_iB_i'X_i+Q_i=0
\end{align}
for an arbitrary index $i\in{\cal N}$.
By the stabilizability of $(A_i,B_i)$,
such an ARE admits the stabilizing solution $X_i\geq 0$,
provided that $Q_i\geq 0$ and $(Q_i^{1/2},A_i)$ is observable for
any unstable modes of $\dot{x}_i(t)=A_ix_i(t)$
corresponding to the imaginary eigenvalues of $A_i$.
The matrix $Q_i$ should be chosen to optimize the
local control performance for each $i$.
It follows that $F_i=B_i'X_i$ is stabilizing, and
more importantly, both
%from \cite{am} (page 106) and \cite{k} (page 238) that
\begin{align} \label{tfis}
    T_{F_i}(s) = F_i(sI-A_i+B_iF_i)^{-1}B_{i},
    \; \; \; \; \;
    T_F(s)={\rm diag}\left\{T_{F_1}(s),\cdots,T_{F_N}(s)\right\}
\end{align}
are strictly positive real (SPR)
\cite{am} (page 106). An application of
Theorem 6.3 in \cite{k} (page 250) concludes
inequality (\ref{prc}), i.e., (\ref{rfd}),
and thus the stability of the feedback MAS.
\pfbox

\begin{remark}\label{remd}
{\em If the knowledge on $D$ is not available, then
the following FI control protocol
\[
     u_i(t) = \delta_{\rm K}(i-i_R)v_i(t)+
    \sum_{j=1}^N a_{ij}\left[v_i(t)-v_j(t)\right]
\]
will have to be used. The above
results in $\dot{x}(t)=\left[A-B{\mathcal M}F\right]x +
B{\mathcal M}r$ for the feedback MAS.
In this case $F_{i\epsilon}=\epsilon^{-1}B_i'X_i$
can be used as an alternative  for each $i\in{\cal N}$
where $X_i\geq 0$ is the stabilizing solution
to (\ref{arexi}). We claim that there exists an $\epsilon\in (0,\ 1)$
such that $(A-B{\cal M}F_{\epsilon})$ is a Hurwitz matrix
where $F_{\epsilon}={\rm diag}(F_{1\epsilon},\cdots,F_{N\epsilon})$.
Indeed inequality (\ref{rfd}) is now replaced by
\[
    \lambda(s):={\rm det}\left(sI-A+
    B{\cal M}F_{\epsilon}\right)\neq 0
    \; \; \;\forall \ {\rm Re}\{s\}\geq 0.
\]
Denote $\underline{D}={\rm diag}(d_1I_{n_1},\cdots,d_NI_{n_N})$.
By the block diagonal form of $A, B,$ and $F_{\epsilon}$, there holds
\begin{align*}
    \lambda(s)= \det\left[\underline{D}\left(sI-A+
    B{\cal M}F_{\epsilon}\right)\underline{D}^{-1}\right]
    ={\rm det}\left(sI-A+B{\cal M}_DD^{-1}F_{\epsilon}\right).
\end{align*}
Since a diagonal $D>0$ exists such that
inequality (\ref{kap}) holds,
$Z={\cal M}_D-I$ satisfies $Z'+Z>0$. The
above inequality is now equivalent to
\begin{align*} %\label{eq:detstab}
    {\rm det}\left(sI-A+BD^{-1}F_{\epsilon}
    +BZD^{-1}F_{\epsilon}\right)\neq 0
    \; \; \; \forall \ {\rm Re}\{s\}\geq 0,
\end{align*}
that is in turn equivalent to
$\det[I+T_{F_{\epsilon}}(s)Z]\neq 0$
$\forall \ {\rm Re}\{s\}\geq 0$ where
\begin{align} \label{tfise}
    T_{F_{\epsilon}}(s)={\rm diag}[T_{F_{1\epsilon}}(s),
    \cdots,T_{F_{N\epsilon}}(s)], \; \; \; \; \;
    T_{F_{i\epsilon}}(s) = d_i^{-1}F_{i\epsilon}(sI-A_i+
    B_id_i^{-1}F_{i\epsilon})^{-1}B_{i}.
\end{align}
Specifically choosing $\epsilon<2/\max\{d_i\}$ ensures that
\[
    A_{F_i}=A_i-B_id_i^{-1}F_{i\epsilon} = A_i -
    B_i(\epsilon d_i)^{-1}B_i'X_{i}
\]
is a Hurwitz matrix in light of the fact that
ARE (\ref{arexi}) can be rewritten as
\[
    A_{F_i}'X_{i} + X_{i}A_{F_i}
    +X_{i}B_i[2(\epsilon d_i)^{-1}-1]B_i'X_{i} + Q_i = 0
\]
and $2(\epsilon d_i)^{-1}-1>0$ for all $i\in{\cal N}$.
The same argument as in the proof of Theorem \ref{th2}
can be used to conclude the SPR of $T_{F_{i\epsilon}}(s)$
and thus $T_{F_{\epsilon}}(s)$ in (\ref{tfis}), implying that
$(A-B{\cal M}F)$ is indeed a Hurwitz matrix.
The above shows that
distributed stabilization can be achieved without
knowledge of $D$ satisfying (\ref{kap}).
However the local optimality is lost
due to the scaling of $\epsilon$ in
$F_{i\epsilon}=\epsilon^{-1}B_i'X_i$,
and it is debatable if the search for $\epsilon>0$
to achieve the stability of $(A-B{\cal M}F_{\epsilon})$
is any simpler than the search for $D$ satisfying (\ref{kap}).
There are strong incentives to search for $D$ rather than $\epsilon$.
\pfbox}
\end{remark}

%\vspace{2mm}

Theorem \ref{th2} provides a sufficient
condition for stabilizability under the distributed
state feedback control. This sufficient condition
becomes necessary for two special cases as shown next.

%\begin{remark}\label{rem2} {\em
\begin{cor}\label{cor0}
Consider state feedback control for the {\rm MAS} over
the directed graph ${\mathcal G}$.
If feedback stability holds for the {\rm MAS}
consisting of either {\rm (i)}\ homogeneous
multi-input unstable agents or {\rm (ii)}
heterogeneous single input unstable agents with $\{A_i\}_{i=1}^N$
having a common unstable eigenvalue,
then the directed graph
${\mathcal G}$ is connected.
\end{cor}

Proof: %Feedback stability implies stabilizability.
For case (i), the homogeneous hypothesis implies
\[
    F(sI-A)^{-1}B=I_N\otimes F_{\rm a}(sI-A_{\rm a})^{-1}B_{\rm a}
\]
where $(A_i,B_i,F_i)=(A_{\rm a},B_{\rm a},F_{\rm a})$
$\forall$ $i\in{\mathcal N}$.
Using the same procedure as in \cite{ldch,om},
the feedback stability condition in (\ref{rfd})
can be shown to be equivalent to
\begin{align*} %\label{hmmasta}
    \det[I+F_{\rm a}(sI-A_{\rm a})^{-1}B_{\rm a}
    \lambda_i({\mathcal M_D})]\neq 0 \; \forall\
    {\rm Re}[s]\geq 0.
\end{align*}
Since $A_{\rm a}$ has unstable eigenvalues by the hypothesis,
the above inequality implies
$\lambda_i({\mathcal M_D})\neq 0$ for all $i$,
concluding that the graph ${\mathcal G}$ is connected.
For case (ii), feedback stability implies
stability of $(A-B{\mathcal M_D}F)$
for some $F$. Thus
\[
    {\rm rank}\left\{\left[
    \begin{array}{ccc}
    sI_n - A && B{\mathcal M_D}
    \end{array}
    \right]\right\} = n \; \forall\ {\rm Re}[s]\geq 0
\]
where $n=n_1+\cdots+n_N$. Recall $A_i$ has dimension
$n_i\times n_i$.
For single input agents, $B{\mathcal M_D}$
has $N$ columns. Taking $s$ to be the common unstable
eigenvalue of $\{A_i\}_{i=1}^N$ implies that
$(sI_n - A)$ has rank $(n-N)$, and thus
$B{\mathcal M_D}$ has rank $N$,
leading to the conclusion of nonsingular ${\mathcal M_D}$ that concludes the proof for (ii).
\pfbox %\end{remark}

\vspace{2mm}

When the states of the MAS are not available for
feedback, a distributed observer can be designed to
estimate the state of each agent, which can then be
used for feedback control. See \cite{ldch,zld} for homogeneous MASs.
We will modify the neighborhood observers
in \cite{zld} for designing distributed output feedback
controllers in the case of heterogeneous MASs
to incorporate the relative information of the
output measurements. Let $\hat{x}_i(t)$ be
estimate of $x_i(t)$, and
\begin{align}\label{exy}
    e_{x_i}(t)=x_i(t)-\hat{x}_i(t), \; \; \; \; \;
    e_{y_i}(t)=y_i(t)-C\hat{x}_i(t),
\end{align}
be the estimation error for $x_i(t)$ and $y_i(t)$,
respectively. There holds
\begin{align}\label{yij0}
    [y_i(t)-y_j(t)] - [\hat{y}_i(t)-\hat{y}_j(t)]
    = e_{y_i}(t)-e_{y_j}(t)=C_ie_{x_i}(t)
    -C_je_{x_j}(t).
\end{align}
The neighborhood observer
modified from \cite{zld} is proposed as follows:
\begin{align}\label{dis}
    \dot{\hat{x}}_i &= A_i\hat{x}_i + B_iu_i %\nonumber
    + d_i\delta_{\rm K}(i-i_R)L_iC_ie_{x_i}
    + d_iL_i\sum_{j=1}^N a_{ij}[C_ie_{x_i}-C_je_{x_j}]  %\\
    %& \; \; + g_i[L_iC_i(\hat{x}_i-x_i)-L_0C_0(\hat{x}_0-x_0)]
    %\; \forall\ i\in{\mathcal N}
\end{align}
for $1\leq i\leq N$. The following result holds.

\begin{theorem}\label{th32}
Suppose that $D$ satisfying $(\ref{kap})$
in {\rm Remark} $\ref{rem1}$
is known, and $(A_i,B_i,C_i)$ is both stabilizable and
detectable for all $i\in{\mathcal N}$.
Then there exist distributed output
feedback stabilizing controllers for the underlying
heterogeneous {\rm MAS}, if the feedback graph is connected.
%For the {\rm MAS} over the time-varying graph,
%the distributed output stabilization requires that
%the graph be uniformly connected with
%adequately small $h>0$ and $(\ref{iirank})$ holds for both $q=p$ and $q=m$.
\end{theorem}

Proof: In light of Theorem \ref{th2},
there exists a FI control protocol
$u(t)={\cal M}_Dr(t)-{\cal M}_DFx(t)$ that stabilizes the
feedback MAS, i.e., $(A-B{\cal M}_DF)$
is a Hurwitz matrix. Since $x(t)$ is unavailable,
a distributed estimator in (\ref{dis}) can be
designed by synthesizing $L_i=Y_iC_i'$ with
$Y_i\geq 0$ being the stabilizing solution to ARE
\begin{align}\label{areyi}
    A_iY_i+Y_iA_i'-Y_iC_i'C_iY_i+\tilde{Q}_i=0
\end{align}
for some $\tilde{Q}_i\geq 0$ that can be used to
optimize the local estimation performance. In addition
the state estimation gain
$L_i$ not only stabilizes the local estimation
error dynamics, but also
satisfies the SPR property for the resulting
\begin{align}\label{tlis}
    T_{L_i}(s)=C_i(sI-A_i+L_iC_i)^{-1}L_i
    \; \; \forall\ i\in {\mathcal N}
\end{align}
which is dual to $T_{F_i}(s)$ in (\ref{tfis})
for the case of state feedback.
Taking difference between the state space
equation $\dot{x}_i=A_ix_i+B_iu_i$
and that in (\ref{dis}) leads to
\begin{align*}
    \dot{e}_{x_i} &= A_ie_{x_i}  - d_i
    L_i\sum_{j=1}^N a_{ij}(C_ie_{x_i}-C_je_{x_j})
    -d_i\delta_{\rm K}(i-i_R)L_iC_ie_{x_i}.
\end{align*}
The above results in the collective error dynamics
described by
\[
    \dot{e}_x(t) = \left[A - L{\cal M}_DC\right]e_x(t),
    \; \; \; \; \; e_x(t)=x(t)-\hat{x}(t),
\]
where ${\cal M}_D$ is the same as that in (\ref{sfclu}).
The proof for Theorem \ref{th2} can be adapted
to show the Hurwitz stability for
$A - L{\cal M}_DC$. Using $\hat{x}(t)$ in place
of $x(t)$ for the FI control protocol results in
\[
    \dot{x}(t)=Ax(t)+B{\cal M}_D[r(t)-F\hat{x}(t)]
    = \left[A - B{\cal M}_DF\right]x(t)+
    B{\cal M}_DFe_x(t)+B{\cal M}_Dr(t).
\]
The overall MAS thus admits the state space description:
\begin{align}\label{ss2}
    \left[\begin{array}{c}
    \dot{x}(t) \\
    \dot{e}_x(t)
    \end{array}\right] &=
    \left[\begin{array}{cc}
    A - B{\mathcal M}_DF & B{\mathcal M}_DF \\
    0 & A - L{\mathcal M}_DC
    \end{array}\right]
    \left[\begin{array}{c}
    x(t) \\
    e_x(t)
    \end{array}\right]+\left[\begin{array}{c}
    B{\mathcal M}_D \\
    0
    \end{array}\right]r(t).
\end{align}
The separation principle for stabilization
holds true as manifested in the
collective dynamics (\ref{ss2}) that concludes the proof.
\pfbox

\vspace{2mm}

If $D$ satisfying inequality (\ref{kap})
in Remark \ref{rem1} is
unknown, then the neighborhood observer
\begin{align}\label{dis1}
    \dot{\hat{x}}_i &= A_i\hat{x}_i + B_iu_i %\nonumber
    + \delta_{\rm K}(i-i_R)L_iC_ie_{x_i}
    + L_i\sum_{j=1}^N a_{ij}[C_ie_{x_i}-C_je_{x_j}]  %\\
    %& \; \; + g_i[L_iC_i(\hat{x}_i-x_i)-L_0C_0(\hat{x}_0-x_0)]
    %\; \forall\ i\in{\mathcal N}
\end{align}
can be employed, resulting
in the same error dynamics described in
(\ref{ss2}), except that ${\cal M}_D$ is replaced by
${\cal M}$. Remark \ref{remd} can be used to
synthesize both the stabilizing state feedback gain
and state estimation gain. The detail is omitted.

In light of Corollary \ref{cor0}, the
sufficient condition in Theorem \ref{th32} for
distributed output stabilizability condition
can be made necessary for SISO heterogeneous
MASs when all agents share some common
unstable poles. We would like to point out that Theorem \ref{th32}
assumes the same communication graph at both input and output.
In practice they can be different from each other which
can give more design freedom. In case that the graphs differ,
then both need to be connected.
Moreover each agent may have different
measurement output from the consensus output.
For simplicity, our paper considers
only the case when the measurement output is the same
as the consensus output; however our design method can be easily adapted to fit to the case when they differ
by replacing $C_i$
%. For this case, the current matrix, $C_i$,
in the state estimator by a different $C_i$
corresponding to the measurement output.

\begin{remark}\label{rmfl} {\bf (Design of control/estimation gains)}\\
{\em If $D$ satisfying inequality (\ref{kap})
in Remark \ref{rem1} is known, the protocols in
(\ref{sfu}) and (\ref{dis}) can be used,
which achieve not only local
optimality but also local robustness. Indeed
many known output feedback controllers,
including the controllers designed using LQG,
${\mathcal H}_{\infty}$ loop shaping,
and LQG/LTR methods, are observer based and satisfy the
required SPR property for
\begin{align}\label{tfli}
    T_{F_i}(s)=F_i(sI-A_i+B_iF_i)^{-1}B_i, \; \; \; \; \;
    T_{L_i}(s)=C_i(sI-A_i+L_iC_i)^{-1}L_i.
\end{align}
While LQG is obvious due to its optimality, the SPR property
is kind of obscure for the other two design methods, which will be
clarified next.

For ${\mathcal H}_{\infty}$ loop shaping based on a
right coprime factorization (a dual
result is presented in \cite{mg}, page 69-72),
$F_i=B_i'X_i$ and $L_i=Y_{i\infty}C_i'$
with $X_i\geq 0$ the stabilizing solution to
the control ARE
\begin{align}\label{xilqr}
    A_i'X_{i}+X_{i}A_i-X_{i}B_iB_i'X_{i}+C_i'C_i = 0,
\end{align}
and $Y_{i\infty}\geq 0$ the stabilizing solution to
the filtering ARE
\begin{align}
    A_iY_{i\infty}+Y_{i\infty}A_i'-
    Y_{i\infty}C_i'C_iY_{i\infty}+B_iB_i'
    +\Delta_i &=0
\end{align}
where $\Delta_i=(\gamma_i^2-1)(I+Y_{i\infty}X_i)B_iB_i'(I+X_iY_{i\infty})$.
Since $\gamma_i>\gamma_{i\rm opt}=\sqrt{1+
\lambda_{\max}(X_iY_i)}\geq 1$,
$\Delta_i\geq 0$ is true. Both
$T_{F_i}(s)$ and $T_{L_i}(s)$
are SPR in light of \cite{am,k}.

For LQG/LTR, $F_i=B_i'X_i$ with
$X_i\geq 0$ the stabilizing solution to
the ARE (\ref{xilqr}), and $L_i=Y_iC_i'$
with $Y_i\geq 0$ the stabilizing solution to the ARE
\[
    A_iY_i+Y_iA_i'-Y_iC_i'C_iY_i+q_i^2\tilde Q_i=0
\]
for some design parameter $q_i>0$ sufficiently large. The matrix
$\tilde Q_i = B_iB_i'$ if $P_i(s)=C_i(sI-A_i)^{-1}B_i$
is minimum phase. Otherwise $\tilde Q_i = B_{im}B_{im}'$ where
$P_i(s)=P_{im}(s)B_{ia}(s)$ with
$P_{im}(s) = C_i(sI-A_i)^{-1}B_{im}$
being the minimum phase part of $P_i(s)$ and $B_{ia}(s)$
satisfying $B_{ia}(-s)'B_{ia}(s)=I$ and containing all
unstable zeros of $P_i(s)$ \cite{zf}. Hence both
$T_{F_i}(s)$ and $T_{L_i}(s)$ in (\ref{tfli}) are
again SPR. However if $D$ in
Remark \ref{rem1} is unknown, and $N$ is too large
to search for it, then the protocols in Remark \ref{remd} for distributed control
and in (\ref{dis1}) for distributed estimation will have to be used, that may destroy the
local optimality and robustness of the feedback MAS.
So there are incentives to search for $D$ whenever
it is possible.
\pfbox}\end{remark}

Remark \ref{rmfl} indicates that the
required SPR property helps to
render not only the local optimality,
but also local robustness for the heterogeneous MAS,
in light of the fact that
LQG/LTR design improves the gain/phase margin
and ${\mathcal H}_{\infty}$ loop shaping
provides the robustness against coprime factor uncertainties.
However it remains unclear for the collective robustness.
The next result is concerned with the existence of static
output stabilizing control law for heterogeneous MASs.

%\vspace*{-1mm}

\begin{cor}\label{co1}
Suppose that $(A_i,B_i,C_i)$ is both stabilizable and
detectable satisfying $\det(C_iB_i)\neq 0$ and
$(A_i,B_i,C_i)$ is strictly minimum phase
for all $i\in{\mathcal N}$. If the feedback graph
exists only at the {\rm MAS} input or output but not at both,
then there exist distributed static output
stabilizing controllers for the underlying
heterogeneous {\rm MAS}, if the feedback graph is
connected.
\end{cor}

%\vspace*{-1mm}

\hspace*{-2mm} Proof: It is shown in
\cite{ggu} that under the hypotheses for
$(A_i,B_i,C_i)$, there exists $K_i$ such that
\[
    u_i(t)=K_iy_i(t)=K_iC_ix_i(t)
\]
is an LQR control law
for the system $\dot{x}_i=A_ix_i+B_iu_i$ with $R=I$.
Since this is true for each $i\in{\mathcal N}$,
the result for distributed state feedback
control can be applied, if the feedback
graph exists at the MAS input. Specifically
$T_F(s)$ in the proof of Theorem \ref{th2}
with $F_i=K_iC_i$ for each $i$ can be made
SPR. The corollary is thus true. \pfbox

\section{Output Consensus}\label{Sec:OC}

This section is focused on output consensus,
assuming that the feedback MAS is stabilized by the
observer-based distributed controllers
developed in the previous
section. Two cases will be investigated.
The first considers the case when
the state vector $x_0(t)$ of the reference model
in (\ref{rm}) is available for feedforward control.
We will show how the feedforward gain $F_{0i}$ can be designed
for $r_i(t)=F_{0i}x_0(t)$ to achieve the output
synchronization as required in (\ref{yi0}) for each
$i\in{\cal N}$. The second assumes $x_0(t)$ unavailable.
We will develop distributed cooperative
observers to estimate the feedforward
control signal $r_i(t)=F_{0i}x_0(t)$ based on relative output
measurements of the neighboring agents, and the
desired output trajectory $y_0(t)=C_0x_0(t)$ that
is transmitted to only one agent.
The following assumption is crucial.

\vspace{2mm}

\noindent {\bf Assumption (a):} \ Agent $i$ has input dimension
$m_i$ no smaller than its output dimension $p$, and
each eigenvalue of $A_0$ is a pole of each column of $P_i(s)$
for all $i\in{\cal N}$.

\vspace{2mm}

The above assumption ensures the right invertibility of
$P_i(s)$ at each eigenvalue of $A_0$. If this assumption
fails, weighting functions $\{W_i(s)\}$ can be used
so that $P_i(s)W_i(s)$ satisfies Assumption (a), and
consensus control design can be carried out for the
weighted agents. Such a method is used widely in
control system design \cite{mg}. The next
result is instrumental.

\begin{lemma}\label{lem3}
Suppose that {\rm Assumption\ (a)} holds,
$(A_i,B_i)$ is stabilizable for all $i$, and the feedback graph is connected.
Then $\{F_i\}_{i=1}^N$ exist such that
$\det[sI-A+B{\cal M}_DF]$ is Hurwitz. Denote
\begin{align}\label{tmts}
    T_{{\cal M}_D}(s) =
    C[sI-A+B{\cal M}_DF]^{-1}B{\cal M}_D, \; \; \; \; \;
    T(s) =C(sI-A+BF)^{-1}B,
\end{align}
and $\Delta(s)=T_{{\cal M}_D}(s)-T(s)$. There holds
\begin{align}\label{skp}
    \lim_{s\rightarrow s_{\kappa}}
    \frac{\Delta(s)}{(s-s_{\kappa})^{\mu_{\kappa}-1}} = 0.
\end{align}
at each eigenvalue $\lambda_{\kappa}$ of
$A_0$ with multiplicity $\mu_{\kappa}$.
\end{lemma}

Proof: Recall $P(s)=C(sI-A)^{-1}B$ and denote
$P_F(s)=F(sI-A)^{-1}B$. The hypotheses imply
\begin{align*}
    \Delta(s)=& C[sI-A+B{\cal M}_DF]^{-1}B{\cal M}_D-C(sI-A+BF)^{-1}B \\
    =& P(s)[I+{\cal M}_DP_F(s)]^{-1}{\cal M}_D
    - P(s)[I+P_F(s)]^{-1} \\
    =& P(s)P_F(s)^{-1} [I + {\cal M}_D^{-1}P_F(s)^{-1}]^{-1}
    -P(s) P_F(s)^{-1}[I+P_F(s)^{-1}]^{-1} \\
    =&P(s)P_F(s)^{-1}\left\{[I + {\cal M}_D^{-1} P_F(s)^{-1}]^{-1}
    -[I+P_F(s)^{-1}]^{-1} \right\}.
\end{align*}
Since $P(s)$ and $P_F(s)$ are block diagonal
transfer matrices and each of their
columns has $s_{\kappa}$ as pole with multiplicity
$\mu_{\kappa}$, $P_F(s)^{-1} \rightarrow 0$,
$[(s-s_{\kappa})^{\mu_{\kappa}-1}P_F(s)]^{-1} \rightarrow 0$,
and $P(s)P_F(s)^{-1}$ approaches a finite block
diagonal matrix as $s \rightarrow s_{\kappa}$. Moreover
\begin{align*}
    [I+P_F(s)^{-1}]^{-1} &= I-P_F(s)^{-1}+o(\{P_F(s)^{-1}\}^2), \\
    [I + {\cal M}_D^{-1}P_F(s)^{-1}]^{-1}=&
    I-{\cal M}_D^{-1}P_F(s)^{-1} +o(\{P_F(s)^{-1}\}^2),
\end{align*}
with $o\left( [P_F(s)^{-1}]^2 \right)$ indicating
that each of its terms approaches zero
in the order of $(s-s_{\kappa})^{2\mu_{\kappa}}$
as $s \rightarrow s_{\kappa}$. Consequently there holds
\begin{align*}
    \Delta(s) \ \rightarrow \ P(s) P_F(s)^{-1}\left\{ I -
    {\cal M}_D^{-1}\right\}P_F(s)^{-1} + o([P_F(s)^{-1}]^2)
\end{align*}
as $s \rightarrow s_{\kappa}$.
Substituting the above into the left hand side of
(\ref{skp}) yields
\begin{align*}
    \frac{\Delta(s)}{(s-s_{\kappa})^{\mu_{\kappa}-1}} \
    \rightarrow \
    \frac{o([P_F(s)^{-1}]^2)}{(s-s_{\kappa})^{\mu_{\kappa}-1}}
    + \frac{P(s) P_F(s)^{-1}\left\{ I -
    {\cal M}_D^{-1}\right\}
    P_F(s)^{-1}}{(s-s_{\kappa})^{\mu_{\kappa}-1}}
    \ \rightarrow \ 0
\end{align*}
as $s\rightarrow s_{\kappa}$, that concludes the proof. \pfbox

\vspace{2mm}

Lemma \ref{lem3} indicates that %under its hypotheses,
$\Delta(s)R(s) = T_{{{\cal M}_D}}(s)R(s)-T(s)R(s)$
has no pole at $s_{\kappa}$. Otherwise its
partial fraction in computing
the term with pole at $s_{\kappa}$
would contradict the limit in (\ref{skp}).
Since $s_{\kappa}$ is an arbitrary eigenvalue of
$A_0$, no eigenvalue of $A_0$ is a pole of $\Delta(s)R(s)$,
implying the stability of $\Delta(s)R(s)$.
The next result is concerned with the
full information (FI) control protocol.

\begin{theorem}\label{th3}
Suppose that {\rm Assumption\ (a)} holds,
the feedback graph is connected, and
the realization $\{A_i,B_i,C_i\}$ is both stabilizable
and detectable for all $i\in{\cal N}$.
Then there exist {\rm FI} control protocols
in the form of $(\ref{sfu})$
such that not only $(A-B{\cal M}_DF)$ is a Hurwitz matrix, but also
the output consensus as required in
$(\ref{yi0})$ is achieved.
\end{theorem}

Proof: Under the assumption of connectivity for the feedback graph and
stabilizability for $(A_i,B_i)$ $\forall\ i$, stabilizing
$\{F_i\}_{i=1}^N$ can be synthesized according
to the proof of Theorem \ref{th2} such that not only
each $(A_i-B_iF_i)$, but also %the global
$A-B{\cal M}_DF$ is a Hurwitz matrix. Hence both
$T(s)$ and $T_{{\cal M}_D}(s)$ as in
(\ref{tmts}) are stable where
\begin{align}\label{ttis}
    T(s)={\rm diag}[T_1(s),\cdots,T_N(s)], \; \; \; \; \;
    T_i(s)=C_i(sI-A_i+B_iF_i)^{-1}B_i.
\end{align}
It is clear that $T_{{\cal M}_D}(s)$ is the transfer matrix from
$r(t)$ to $y(t)$ for the feedback MAS under the
FI control protocol, and it becomes equal to
$T(s)$, if ${\cal M}_D=I$, corresponding to
$N$ decoupled feedback systems under the FI control.
In addition the hypotheses on stabilizable and detectable
realization for $P_i(s)$ and on eigenvalues
of $A_0$ being poles of each column of $P_i(s)$
imply that none of the eigenvalues of
$A_0$ is a transmission zero of $P_i(s)$
for all $i$. Consequently the two equations
\begin{align}\label{sylv}
    %\begin{array}{rcl}
    (A_i-B_iF_i)\Pi_i + B_iF_{0i} = \Pi_iA_0, \; \; \; \; \;
    C_i\Pi_i - C_{0} = 0,
    %\end{array}
\end{align}
admit solutions $(\Pi_i,F_{0i})$ for each $i$, in light of
the internal model principle \cite{h}. Let $R_i(s)$ and $Y_0(s)$
be the Laplace transform of $r_i(t)=F_{0i}x_0(t)$
and $y_0(t)=C_0x_0(t)$, respectively. It is easy to
see that $T_i(s)R_i(s)-Y_0(s)$ is stable,
and the output of $T_i(s)$ tracks $y_0(t)=C_0x_0(t)$
with zero steady-state error for all $i$.
To show the output consensus for the
feedback MAS under the FI control protocol,
we examine the synchronization error.
Recall $\Delta(s)=T_{{\cal M}_D}(s)-T(s)$. Denote
$\underline{Y}_0(s)$ and $R(s)$ as the Laplace transform of
$\underline{y}_0(t)={\bf 1}_N\otimes y_0(t)$,
and $r(t)=F_0\underline{x}_0$, respectively
with $\underline{x}_0(t)={\bf 1}_N\otimes x_0(t)$.
Then the synchronization error in the $s$-domain is given by
\begin{align}\label{eys}
    E_y(s)=T_{{\cal M}_D}(s)R(s)-\underline{Y}_0(s)=
    \Delta(s)R(s)+[T(s)R(s)-\underline{Y}_0(s)]
\end{align}
%is the synchronization error in the $s$-domain
for the feedback MAS under the FI control protocol (\ref{sfclu}).
Since the first term on the right hand side is stable
by Lemma \ref{lem3}, and the
second term on right of (\ref{eys}) is also
stable, $e_y(t)=y(t)-\underline{y}_0(t)
\rightarrow 0$ as $t\rightarrow\infty$, thereby concluding the proof.
\pfbox

\begin{remark}\label{obr}{\em
If the state vector of each agent is not
available for feedback control,
then the state estimation results in
the previous section can be employed.
Specifically the estimator in (\ref{dis})
can be employed to result
in the error dynamics described by
\begin{align*}
    \left[\begin{array}{c}
    \dot{x}(t) \\
    \dot{e}_x(t)
    \end{array}\right] =
    \left[\begin{array}{cc}
    A - B{\mathcal M}_DF &  B{\mathcal M}_DF\\
    0 & A - L{\mathcal M}_DC
    \end{array}\right]
    \left[\begin{array}{c}
    x(t) \\
    e_x(t)
    \end{array}\right]+
    \left[\begin{array}{c}
    B{\mathcal M}_D \\
    0
    \end{array}\right]r(t)
\end{align*}
The output consensus can be achieved under
the same hypotheses as those in Theorem \ref{th3},
in light of the fact that
$e_x(t)\rightarrow 0$ as $t\rightarrow\infty$,
provided that $(A_i-L_iC_i)$ is Hurwitz for all
$i\in{\cal N}$. If the state vector of the reference
model is not available for feedforward control, then
the estimator design is more involved, which will
be tackled in the later part of the section.
\pfbox}\end{remark}

Assumption (a) is not only important
for achieving output consensus, but also has
an important implication that is valuable
to output estimation in the case when
full information is not available
for output consensus. This is presented in the next result.

\begin{cor}\label{cor1}
Suppose that $T_{F_i}(s)$ in $(\ref{tfis})$
is not only internally stable but also {\rm SPR} for all $i\in{\cal N}$,
$M_D$ satisfies inequality $(\ref{LE1})$, and the rest of the
hypothesis of {\rm Theorem} $\ref{th3}$ holds. Then the feedforward
gains $\{F_{0i}\}_{i=1}^N$ achieving the
output consensus satisfy %such that
$A_i\Pi_i=\Pi_iA_0$ and $F_{0i}=F_i\Pi_i$
for all $i\in{\cal N}$.
\end{cor}

Proof: In light of Theorem \ref{th3}, there
exists a FI control protocol that achieves not
only feedback stability but also output
consensus. Hence there exists a solution pair
$(\Pi_i,F_{0i})$ to (\ref{sylv}) for each $i\in{\cal N}$. Denote
$\underline{C}_0=I_N\otimes C_0$,
$\underline{A}_0=I_N\otimes A_0$, and
$\Pi={\rm diag}(\Pi_1,\cdots,\Pi_N)$.
Then (\ref{sylv}) for all
$i\in{\cal N}$ can be put together, leading to
\begin{align}\label{exch}
    \begin{array}{rl}
    (A-BF)\Pi+BF_0 =\Pi\underline{A}_0, \; \; \; \; \;
    C\Pi -\underline{C}_0 =0.
    \end{array}
\end{align}
Recall $E_y(s)$ in (\ref{eys}) under the FI control protocol
as in the proof of Theorem \ref{th3}. Its
dynamic behaviors in the time domain are described by
\begin{align}\nonumber
    \left[\begin{array}{c}
    \dot{x}(t) \\ \dot{\underline{x}}_0(t)
    \end{array}\right]  &=
    \left[\begin{array}{cc}
    A-B{\cal M}_DF & B{\cal M}_DF_0 \\ \label{tmfs}
    0 & \underline{A}_0
    \end{array}\right]\left[\begin{array}{c}
    x(t)\\ \underline{x}_0(t)
    \end{array}\right], \\
    e_y(t) &= Cx(t) - \underline{C}_0\underline{x}_0(t)
\end{align}
where $\underline{x}_0(t)={\bf 1}_N\otimes x_0(t)$,
and $e_y(t)$ and $E_y(s)$ are a Laplace transform pair.
In light of
Theorem \ref{th3}, we can show that $e_y(t)%=C\tilde{x}(t)
\rightarrow 0$ as $t\rightarrow\infty$.
By the hypothesis on $M_D$,
${\cal M}_D=I+Z$ and $Z+Z'>0$.
%and $Z_q$ is a positive definite operator.
The above state equation can now be written as
%Substituting ${\cal M}_q=I+Z_q$ into the above
%state equation yields
\begin{align}\label{xx0}
    \left[\begin{array}{c}
    \dot{x}(t) \\ \underline{\dot{x}}_0(t)
    \end{array}\right]  &=
    \left[\begin{array}{cc}
    A-BF & BF_0 \\
    0 & \underline{A}_0
    \end{array}\right]\left[\begin{array}{c}
    x(t)\\ \underline{x}_0(t)
    \end{array}\right]
    +\left[\begin{array}{c}
    BZ\\ 0
    \end{array}\right][r(t)-Fx(t)].
\end{align}
Applying the similarity transform
\[
    \left[\begin{array}{c}
    \tilde{x}(t)\\ \underline{x}_0(t)
    \end{array}\right]
    = \left[\begin{array}{cr}
    I & -\Pi \\ 0 & I
    \end{array}\right]
    \left[\begin{array}{c}
    x(t)\\ \underline{x}_0(t)
    \end{array}\right]
\]
to state equation (\ref{xx0}) and error
equation (\ref{tmfs}), and utilizing the relation
$C\Pi -\underline{C}_0=0$ yield
\begin{align*} \nonumber
    \dot{\tilde{x}}(t) =(A-BF)\tilde{x}(t)+
    BZ[r(t)-Fx(t)], \; \; \; \; \; \;
    \underline{\dot{x}}_0(t)=\underline{A}_0
    \underline{x}_0(t),
\end{align*}
and $e_y(t) = C\tilde{x}(t)$.
By $r(t)=F_0\underline{x}_0(t)$, there holds
\begin{align}\label{vtt}
    v(t) :=\ r(t)-Fx(t)=F_0\underline{x}_0(t) - Fx(t)
    = (F_0-F\Pi)\underline{x}_0(t)-F\tilde{x}(t).
\end{align}
As a result, the state equation for
$\tilde{x}(t)$ can be written as
\[
    \dot{\tilde{x}}(t)= (A-B{\cal M}_DF)\tilde{x}(t)+
    BZ(F_0-F\Pi)\underline{x}_0(t).
\]
%The transfer matrix from input
%$Z(F_0-F\Pi)\underline{x}_0(t)$ to output $e_y(t)$
%is $C(sI-A+B{\cal M}_DF)^{-1}B$ that
%is internally stable, and
Since $e_y(t)=C\tilde{x}(t)\rightarrow 0$
as $t\rightarrow\infty$ and  $(A-B\mathcal M_DF)$ is internally stable, %, in light of
%Theorem \ref{th3}.
it holds that
\[
    (F_0-F\Pi)\underline{x}_0(t) \ \rightarrow \ 0
\]
as $t\rightarrow\infty$ in light of %the strict positivity of
$Z+Z'>0$. The persistency or divergence of
$\underline{x}_0(t)$ imply that
$(F_0-F\Pi)=0$, concluding $F_{0i}=F_i\Pi_i$
$\forall i\in{\cal N}$. Upon substituting $F_0=F\Pi$
into (\ref{exch}) yields $A\Pi = \Pi\underline{A}_0$
and thus $A_i\Pi_i=\Pi_iA_0$ for all $i\in{\cal N}$,
that concludes the proof. \pfbox

%\vspace{2mm}

\begin{remark}{\em The solution $\Pi_i$
to $A_i\Pi_i=\Pi_iA_0$ and $F_{0i}=F_i\Pi_i$
for each $i$ does not have to be solved from
(\ref{exch}). Once $F_{0i}$ and $F_i$ are available,
$\Pi_i$ can be computed via
\[
    \left[\begin{array}{c}
    C_i\\ F_{i}
    \end{array}\right]\Pi_i = \left[\begin{array}{c}
    C_0\\ F_{0i}
    \end{array}\right] \; \; \;
    \Longrightarrow \; \; \; \Pi_i =
    \left[\begin{array}{cc}
    C_i' & F_{i}'
    \end{array}\right]
    \left(\left[\begin{array}{c}
    C_i\\ F_{i}
    \end{array}\right]
    \left[\begin{array}{cc}
    C_i' & F_{i}'
    \end{array}\right]\right)^{-1}
    \left[\begin{array}{c}
    C_0\\ F_{0i}
    \end{array}\right],
\]
provided that the inverse exists. If the inverse does
not exist, pseudo-inverse can be used in place of
the inverse. It is also important to point out that
$\{F_{0i}\}$ can be synthesized rather easily which
will be demonstrated in the example section. }\pfbox
\end{remark}

In practice it is possible that
the states of the MAS are not
available for feedback control and the state
of the reference model is not available
for feed-forward control, giving rise to
the problem of distributed observer.
By Corollary \ref{cor1} and (\ref{vtt}), there holds
\begin{align*}
    v(t)&=-F\tilde{x}(t), \; \; \;
    \tilde{x}(t)=x(t)-\Pi\underline{x}_0(t).
\end{align*}
Consequently the FI control protocol (\ref{sfu})
for $1\leq i\leq N$ yields %is equivalent to
\begin{align} \label{utte}
    u(t) = {\cal M}_Dv(t) %[F_0\underline{x}_0(t)-Fx(t)]
    =-{\cal M}_DF\tilde{x}(t).
\end{align}
In addition the result on $A\Pi=\Pi\underline{A}_0$
in Corollary \ref{cor1} shows that
%in the proof of Corollary \ref{cor1} implies that
\begin{align}\label{xtte}
    %\begin{array}{rl}
    \dot{\tilde{x}}(t)=A\tilde{x}(t)+Bu(t), \; \; \;
    e_y(t)=C\tilde{x}(t),
    %\end{array}
\end{align}
by $\dot{\tilde{x}}(t)=Ax(t)
+Bu(t)-\Pi\underline{A}_0\underline{x}_0
= Ax(t)+Bu(t)-A\Pi\underline{x}_0$.
Hence we have an estimation problem
for $\tilde{x}(t)$, governed by
the state space system in (\ref{xtte}),
rather than for both $x(t)$ and $x_0(t)$.
As such the relative output measurement
for agent $i$ in (\ref{yij0}) can be expressed as
\begin{align}\label{eyij}
    \varepsilon_{i,j}(t)=y_i(t)-y_j(t)
    = e_{y_i}(t)-e_{y_j}(t)=C_i\tilde{x}_i(t)-C_j\tilde{x}_j(t).
\end{align}
%for those $j\in{\cal N}$ being neighbors of the
%$i$th agent.
Recall that $e_{y_i}(t)=y_i(t)-y_0(t)$ is the
tracking error for the $i$th agent, i.e., the $i$th
component of $e_y(t)$. There is at least
one agent that has the access to
\begin{align}\label{eyir}
    e_{y_{i_R}}=y_{i_R}(t)-C_0x_0(t)=C_{i_R}x_{i_R}(t)-C_0x_0(t)
    = C_{i_R}\tilde{x}_{i_R}(t)
\end{align}
with index $i_R$ corresponding to a reachable node.
The above leads to the distributed estimator:
\begin{align}\label{dist}
    \dot{\hat{\tilde{x}}}_i =
    A_i\hat{\tilde{x}}_i + B_iu_i %\nonumber
    + d_i\delta_{\rm K}(i-i_R)L_iC_ie_{\tilde{x}_i}
    + d_iL_i\sum_{j=1}^N a_{ij}
    [C_ie_{\tilde{x}_i}-C_je_{\tilde{x}_j}]  %\\
    %& \; \; + g_i[L_iC_i(\hat{x}_i-x_i)-L_0C_0(\hat{x}_0-x_0)]
    %\; \forall\ i\in{\mathcal N}
\end{align}
modified from (\ref{dis}) with
$e_{\tilde{x}_k}(t)=\tilde{x}_k(t)-\hat{\tilde{x}}_k(t)$.
The next result is parallel to Theorem \ref{th32}.

\begin{theorem}\label{thf}
Suppose that the feedback graph ${\cal G}$
is connected, and each agent
admits a stabilizable and detectable realization.
Let $D$ satisfying $(\ref{kap})$ be known.
Then there exist $F_i$ and $L_i$ such that
both %$(A_i-B_iF_i)$ and $(A_i-L_iC_i)$
$T_{F_i}(s)$ and $T_{L_i}(s)$ in $(\ref{tfli})$ are
{\rm SPR}, and %are Hurwitz matrix. In addition,
the distributed estimator $(\ref{dist})$ and protocol
$u(t)=-\mathcal M_DF\hat{\tilde x}(t)$
achieve output consensus.
\end{theorem}

Proof: Taking difference between $\dot{x}_i(t)=A_ix_i(t)+B_iu_i(t)$
and local estimator (\ref{dist}) yields
\[
    \dot{e}_{\tilde{x}_i} =
    A_ie_{\tilde{x}_i}
    + d_i\delta_{\rm K}(i-i_R)L_iC_ie_{\tilde{x}_i}
    + d_iL_i\sum_{j=1}^N a_{ij}
    [C_ie_{\tilde{x}_i}-C_je_{\tilde{x}_j}]
\]
for each $i\in{\cal N}$. Packing
together leads to the collective error equation of
\[
    \dot{e}_x(t) = (A-L{\cal M}_DC)e_x(t).
\]
It is now easy to verify that with
$u(t)=-\mathcal M_DF\hat{\tilde x}(t)$, the overall
feedback MAS admits the state space equation
\begin{align}\label{ss3}
    \left[\begin{array}{c}
    \dot{\tilde{x}}(t) \\
    \dot{e}_{\tilde{x}}(t)
    \end{array}\right] &=
    \left[\begin{array}{cc}
    A - B{\mathcal M}_DF & B{\mathcal M}_DF \\
    0 & A - L{\mathcal M}_DC
    \end{array}\right]
    \left[\begin{array}{c}
    \tilde{x}(t) \\
    e_{\tilde{x}}(t)
    \end{array}\right].
\end{align}
In light of the hypothesis on %the SPR
$T_{F_i}(s)$ and $T_{L_i}(s)$ being SPR,
%as given in $(\ref{tfli})$,
both $A - B{\mathcal M}_DF$ and
$A - L{\mathcal M}_DC$ are Hurwitz.
Recall the proof of Theorem \ref{th2}.
It follows that $\tilde{x}(t)\rightarrow0$ and $e_{\tilde{x}}(t)\rightarrow0$ as $t\rightarrow\infty$. Consequently
$\tilde x(t)\rightarrow 0$ as $t\rightarrow\infty$.
In light of the relations
$C_0=C_i\Pi_i$ in (\ref{exch}) and
$\tilde x_i(t)=x_i(t)-\Pi_ix_0(t)$, there holds
\[
    C_i\tilde x_i(t)=
    C_i[x_i(t)-\Pi_ix_0(t)] =
    C_ix_i(t)-C_0x_0(t)=
    y_i(t)-y_0(t) \ \rightarrow \ 0
\]
as $t\rightarrow\infty$. Hence the output consensus is achieved.
\pfbox

%\vspace{2mm}

\begin{remark}{\em It is important to note
that due to equation (\ref{eyir}),
the reference trajectory $y_0(t)=C_0x_0(t)$ is
transmitted to only one node, i.e., one of the reachable
nodes with index $i_R$. This has the advantage of lowering the
communication overhead. However, there is a tradeoff between
consensus robustness and communication
overhead. If $y_0(t)$ is transmitted to
more than one nodes, then the single point of
failure scenario is avoided. This way
the closed loop system is more robust in the
presence of broken communication links.
There exists a tradeoff between the communication
overhead and consensus robustness. }\pfbox
\end{remark}

\section{An Illustrative Example}
\begin{figure}[!htb]
\centering
	\includegraphics[scale=.7]{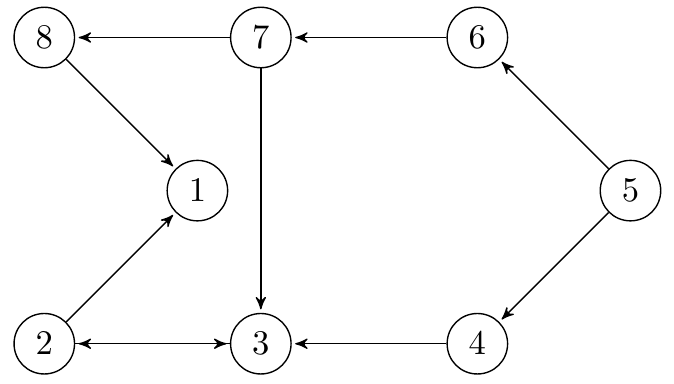}
\caption{Graph for $N=8$ point masses.}
	\label{Fig:Graph}
\end{figure}

In this section we
will demonstrate with a numerical example the output consensus results
presented in this paper.
Following \cite{ls}, consider a system of $8$ point masses
moving in one spatial dimension.
Dynamics are governed by
\begin{align*}
    \dot x_i&=A_ix_i+B_iu_i, \; \; \; \; \;
    y_i=C_ix_i
\end{align*}
where $A_i=\begin{bmatrix}0&1\\0&
    -f_{d_i}\end{bmatrix}$, $B_i=\begin{bmatrix}0\\ \frac{1}{m_i}\end{bmatrix}$, and $C_i=\begin{bmatrix}1&0\end{bmatrix}$ for $i=1,2,3,7,8$; and
\begin{align*}
    \dot x_i&=\frac{1}{m_i}u_i, \; \; \; \; \;
    y_i=x_i
\end{align*}
for $i=4,5,6$. The first set of dynamics represents agents that
experience drag forces and whose acceleration
is directly controlled, while the second set represents
agents whose velocity is directly controlled. The output
signal corresponds to the position of the point mass.
Figure \ref{Fig:Graph} shows the interconnection graph
for the network of 8 agents, with parameters
$\left\{m_i\right\}=\left\{2.5,3,4,1,2,5,7,6\right\}$ and
$\left\{f_{d_i}\right\}=\left\{0.5,0.9,1.9,0,0,0,1.1,3.9\right\}$.

\begin{figure}[!htb]
\centering
	\includegraphics[scale=.74]{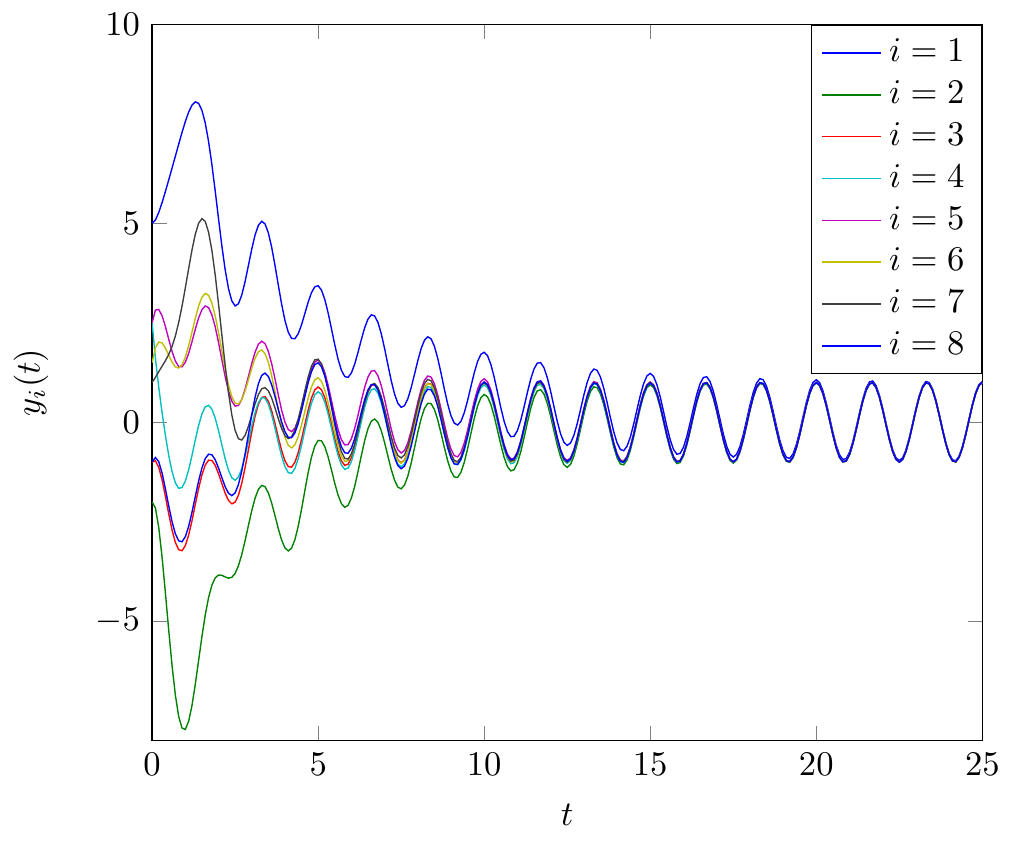}
\caption{Output consensus using a full information control law.}
	\label{Fig:FIsim}
\end{figure}

\begin{figure}[!htb]
\centering
	\includegraphics[scale=.75]{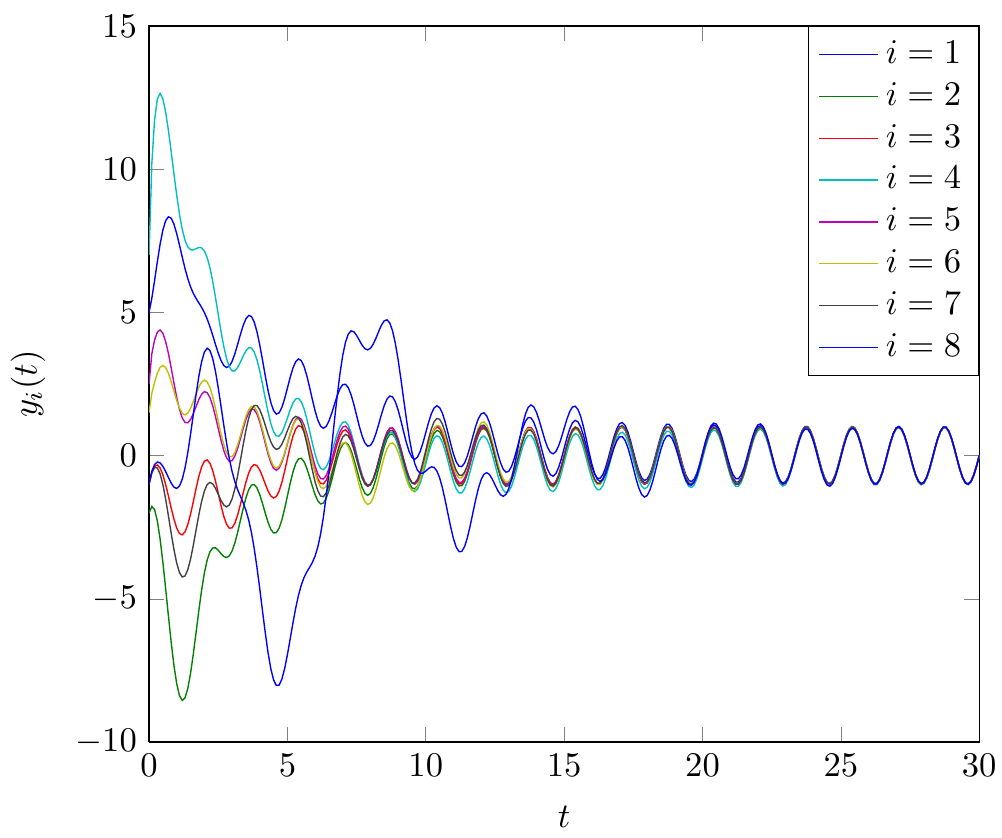}
\caption{Output consensus using a full information control law with distributed estimator.}
	\label{Fig:FIobsvsim}
\end{figure}

It is important to note that the network graph is connected
with only one reachable node at node 1.
Consequently condition (\ref{LE1})
is satisfied with $i_R=1$. This means that only
agent 1 has direct access to its own full information control signal,
every other agent only has access to the difference between its and
its neighbor's full information control signal. The consensus problem
for this example is to synchronize all agents to the signal
$y_0(t)=\cos(\omega_0t)$ generated by the reference model (\ref{rm})
with $\omega_0=3.77$. To satisfy assumption (a), we use a weighting function, $W(s)$, with poles at $\pm j\omega_0$ so that $P_i(s)W(s)$
(where $P_i(s)$ is the transfer function of agent $i$)
contains the eigenvalues of the reference model, i.e.,
$\pm j\omega_0$. In addition, $W(s)$ includes
a proportional-derivative (PD) compensator
by adding a finite stable zero to its dynamics
which helps to enlarge the bandwidth and thus
enhance the transient performance of the control system.
The state feedback and estimation gains are computed for each
agent individually by solving ARE (\ref{xilqr}) and its dual,
respectively. The feedforward gain for each agent could be computed
by solving (\ref{sylv}) for each agent. However, we take this opportunity to present a simple method that may be used in place of solving (\ref{sylv}): The feedfoward gain has the form $F_{0_i}=\begin{bmatrix}F_{0_{i_1}}&F_{0_{i_2}}\end{bmatrix}$. Denoting $X_0(s)$ as the Laplace transform of $x_0(t)$, it is easy
to see that for each agent in closed loop,
\begin{align*}
    Y_i(s)=T_i(s)F_{0_i}X_0(s)&=
    T_i(s)\left(\frac{\omega_0F_{0_{i_1}}}{s^2+\omega_0^2}
    +\frac{sF_{0_{i_2}}}{s^2+\omega_0^2}\right)
    =\frac{K}{s-j\omega_0}+\frac{K^\ast}{s+j\omega_0}
\end{align*}
by partial fraction expansion. Thus, to enforce $y_i(t)=\cos(\omega_0t)$ and achieve tracking we must choose $K$ such that $|K|=1$ and $\angle K=0$. In other words, $F_{0_i}$ must be chosen such that
\begin{align*}
    \sqrt{F_{0_{i_1}}^2+F_{0_{i_2}}^2}&=\frac{1}{T_i(j\omega_0)},
    \; \; \; \; \; \angle  \left[F_{0_{i_1}}
    +jF_{0_{i_2}}\right]=\frac{\pi}{2}-\angle T_i(j\omega_0)
\end{align*}
admit a solution $F_{0_{i_1}}$ and $F_{0_{i_2}}$. Indeed, it is always possible to find a unique solution to both equations.
It is important to remember
that control design, either using our proposed method or by solving (\ref{sylv}), can be carried out in a fully distributed manner. Furthermore, the simple method we have proposed to compute $F_{0_i}$ is applicable to tracking step and ramp functions, and sinusoids with arbitrary amplitudes and phase angles. In fact, it is applicable when the signal to track consists of any linear combination of the previously mentioned functions.
Figure \ref{Fig:FIsim} shows that output consensus is achieved when
the state vector of each agent and the reference model is available.
For the case when the distributed estimator in (\ref{dist}) is required, Figure
\ref{Fig:FIobsvsim} illustrates that consensus is also achieved.

\begin{figure}[!htb]
\centering
	\includegraphics[scale=.75]{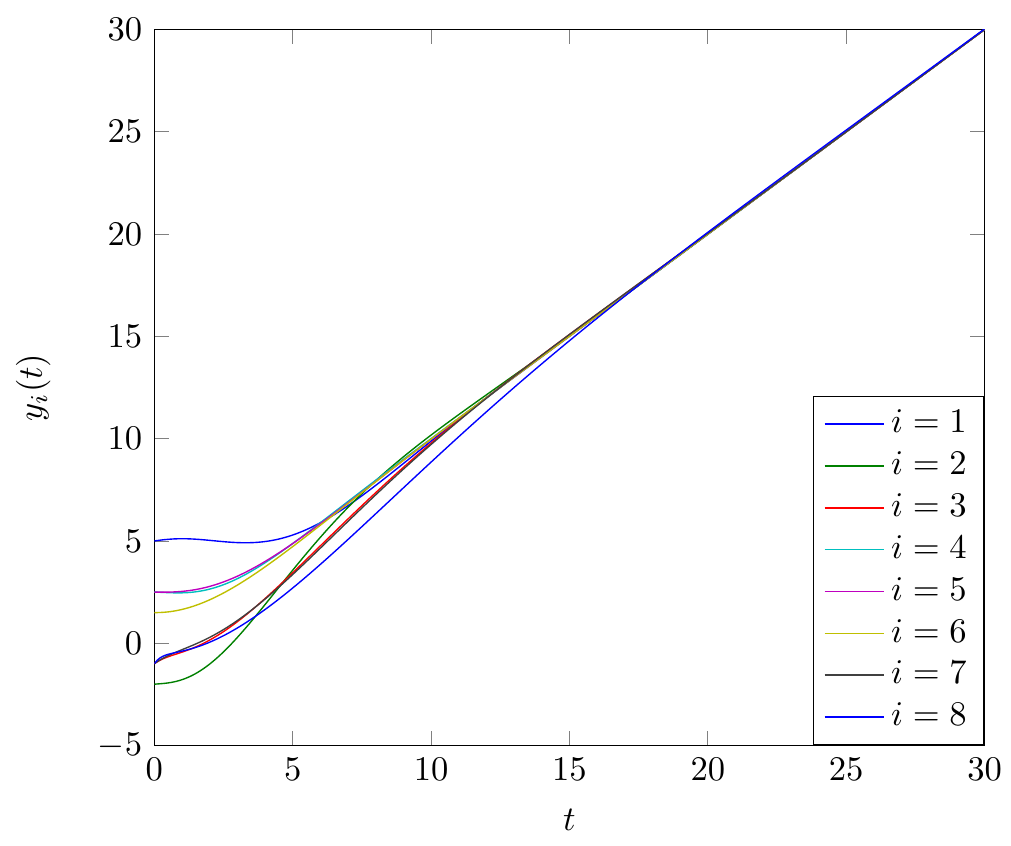}
\caption{Output consensus using a full information control law.}
	\label{Fig:FIsimramp}
\end{figure}

\begin{figure}[!htb]
\centering
	\includegraphics[scale=.75]{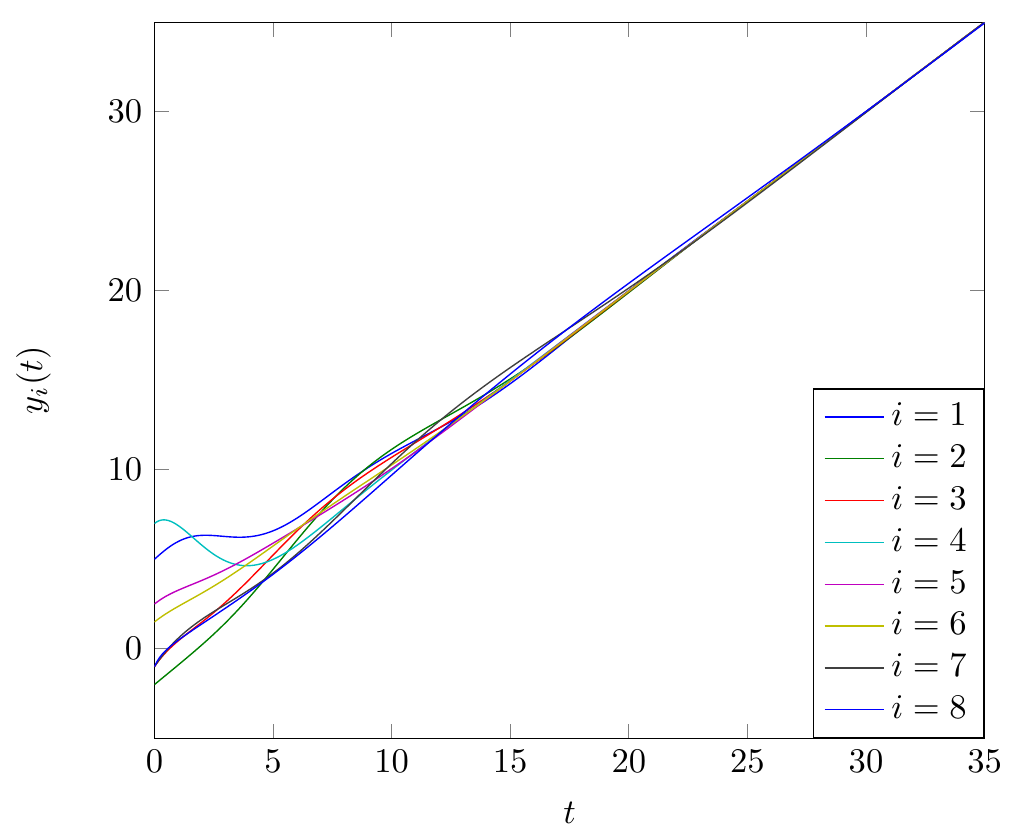}
\caption{Output consensus using a full information control law with distributed estimator.}
	\label{Fig:FIobsvsimramp}
\end{figure}

The simulation study is also carried out
for the unit ramp reference signal. 
Figure \ref{Fig:FIsimramp} shows that output
consensus is achieved under the FI control protocol when the state vector
of each agent and the reference model is available. 
When the states of the $N$ agents and the reference model
are unavailable, the FI control signal can be estimated using the 
distributed estimator in (\ref{dist}). Figure \ref{Fig:FIobsvsimramp}
shows that output consensus is achieved again.
%However the difference is minor. 

For both sinusoidal and ramp reference signals, the consensus
control based on the estimator in (\ref{dist}) takes a
longer time to synchronize the output signals. In addition, there are greater 
fluctuations in the initial stage of the
consensus process than that based on the FI control protocol. 
The reason lies in the fact
that the estimation error affects 
the tracking performance negatively, and full synchronization
is not possible until the the estimation error is 
negligibly small. Certainly, there are incentives to design
faster estimators for consensus control.

\section{Conclusion}

Output consensus control for continuous-time
heterogeneous MASs is studied in this paper,
aimed at synchronizing all the agents' output
to the desired trajectory generated by a reference model.
Our contribution includes the use of the positive real
real property of transfer matrices in achieving the distributed stability
for the feedback MAS, and the
synthesis algorithm developed for design of
the consensus protocol to achieve output consensus.
We have demonstrated that output
consensus control can be decomposed into
tracking regulation for each
individual agent control system, provided that
the graph is connected and each agent's
dynamics contain the reference dynamics
as its internal modes. In establishing our
main results, Lemma 1 and the
positivity of the graph matrix ${\cal M}_D$
in Section 2 play a fundamental role.
In order to achieve consensus to a reference trajectory,
it is sufficient for one agent to have access
to the reference signal, which lowers the communication
overhead for the MAS. %, however, it is advised that several
%agents have access to the reference to avoid the single
%point of failure problem.
In addition it is not necessary to
duplicate the reference model in each of the
$N$ local and distributed feedback controllers, thereby
eliminating synchronization of the local reference
models commonly required in the existing work
for consensus control. Thus the
communication cost can be lowered further. Furthermore the
communication graph can be different at the input and output.
Our controller synthesis is based on $\mathcal
H_\infty$ loop shaping and
LQG/LTR methods, and therefore can accommodate
performance and robustness requirements.
%Finally, each controlled agent is robust
%to perturbations in the form of coprime factor
%uncertainties ($\mathcal H_\infty$ loop shaping)
%or gain/phase uncertainties (LQG/LTR).
How these local properties translate to the
robustness and performance of the collective
dynamics is currently under study.

\end{document}